\input amstex\documentstyle{amsppt}  
\pagewidth{12.5cm}\pageheight{19cm}\magnification\magstep1
\topmatter
\title Almost special representations of Weyl groups
\endtitle
\author G. Lusztig\endauthor
\address{Department of Mathematics, M.I.T., Cambridge, MA 02139}\endaddress
\thanks{Supported by NSF grant DMS-2153741}\endthanks
\endtopmatter   
\document
\define\mat{\pmatrix}
\define\endmat{\endpmatrix}

\define\mpb{\medpagebreak}

\define\hW{\hat W}

\define\si{\sim}

\define\sqc{\sqcup}

\define\bX{\bar X}

\define\lb{\linebreak}

\define\op{\oplus}
   
\define\part{\partial}
\define\emp{\emptyset}

\define\m{\mapsto}
\define\do{\dots}

\define\sub{\subset}    

\define\T{\times}

\define\nl{\newline}

\define\un{\underline}

\define\Ind{\text{\rm Ind}}

\redefine\spa{\spadesuit}

\define\a{\alpha}
\redefine\b{\beta}

\define\g{\gamma}
\redefine\d{\delta}
\define\e{\epsilon}

\define\r{\rho}
\define\s{\sigma}

\define\k{\kappa}
\redefine\l{\lambda}

\redefine\G{\Gamma}
\redefine\D{\Delta}

\define\Si{\Sigma}

\redefine\L{\Lambda}
\define\Ph{\Phi}

\redefine\ss{\bold s}

\define\CC{\bold C}

\define\NN{\bold N}

\define\ZZ{\bold Z}

\define\ca{\Cal A}

\define\cc{\Cal C}

\define\ci{\Cal I}

\define\cl{\Cal L}

\define\cz{\Cal Z}

\head Introduction\endhead
\subhead 0.1\endsubhead
Let $W$ be an irreducible Weyl group and let $\hW$ be the
set of (isomorphism classes of) irreducible representations
of $W$ over $\CC$. Let $c\sub\hW$ be a family of $W$ (see
\cite{L79},\cite{L84}). Let $\G_c$ be the finite group
associated in \cite{L84} to $c$ and let
$c@>>>M(\G_c)$ be the imbedding defined in \cite{L84}.
(For any finite group $\G$, we denote by $M(\G)$ the set of
$\G$-conjugacy of pairs $(x,\s)$ where $x\in\G$ and $\s$ is
an irreducible representation over $\CC$ of the centralizer
of $x$ in $\G$.)

In \cite{L82} a set $Con_c$ of
(not necessarily irreducible) representations of $W$ with
all irreducible components in $c$ was defined and it was
conjectured that these are exactly the representations of $W$
carried by the various left cells \cite{KL} contained in
the two-sided cell associated to $c$. (This conjecture was
proved in \cite{L86}.)

We can view $Con_c$ as a subset of the $\CC$-vector space
$\CC[c]$ (with basis $c$) and hence with
a subset of the $\CC$-vector space $\CC[M(\G_c)]$ (with
basis $M(\G_c)$) via the imbedding $\CC[c]\sub\CC[M(\G_c)]$
induced by $c\sub M(\G_c)$.
Note that $\Ind_1^{\G_c}(1)$ can be also viewed as an
element of $\CC[M(\G_c)]$ (namely $\sum_\r\dim(\r)(1,\r)$
where $\r$ runs over the irreducible representations of
$\G_c$). We have $\Ind_1^{\G_c}(1)\in Con_c$ except when

(a) $|c|$ equals $2,4,11$ or $17$;
\nl
(in these cases $W$ is of exceptional type). To remedy this,
we enlarge $Con_c$ to the subset
$Con^+_c=Con_c\cup\Ind_1^{\G_c}(1)$ of $\CC[M(\G_c)]$.
(We have $Con^+_c=Con_c$ whenever $c$ is not as in (a)).

The main result of this paper is a definition of

{\it a subset $\ca_{\G_c}\sub c$ in canonical bijection
with $Con^+_c$ such that each element of $\ca_{\G_c}$ 
appears with nonzero coefficient in the corresponding
element of $Con^+_c$.}

In \cite{L79a} a specific
representation $sp_c\in c$ in $c$ was defined (it was later
called the {\it special representation}); it corresponds
to $(1,1)\in M(G_c)$. One of its
properties is that it appears with coefficient $1$ in any
element of $Con^+_c$. We always have $sp_c\in\ca_{\G_c}$. In
fact, $sp_c$ corresponds to $\Ind_1^{\G_c}(1)\in Con^+_c$.
Thus the representations in $\ca_{\G_c}$ generalize $sp_c$;
we call them {\it almost special} representations of $W$.
(This name is justified in 2.4.)

We will show elsewhere that (in the case where $W$ is of
simply laced type) the irreducible representation
of $W$ attached in \cite{L15} to a stratum of $G$ is
almost special.

\subhead 0.2\endsubhead
Our definition of $\ca_{\G_c}$ relies on the theory of
new basis \cite{L19},\cite{L20},\cite{L23}. 

Let $\cz_{\G_c}$ be the set of
pairs $(\G'\sub\G'')$ of subgroups of $\G_c$ with $\G'$
normal in $\G''$. For $(\G'\sub\G'')\in\cz_{\G_c}$ let
$$\ss_{\G',\G''}:\CC[M(\G''/\G')]@>>>\CC[M(\G_c)]$$
be the $\CC$-linear map defined in \cite{L20, 3.1}.
In \cite{L23, 2.3} to $c$ we have associated a subset
$X_{\G_c}$ of $\cz_{\G_c}$.

Let $\bX_{\G_c}=X_{\G_c}\cup\{(S_1,S_1)\}$.
(We denote by $S_n$ the symmetric group in $n$ letters.)
We have $(S_1,S_1)\in X_{\G_c}$ if and only if
$c$ is not as in 0.1(a); in these cases we have 
$\bX_{\G_c}=X_{\G_c}$.
Now,

(a) $(\G'\sub\G'')\m\ss_{\G',\G''}(1,1)$
{\it is a bijection of $\bX_{\G_c}$ onto a subset 
$S(\G_c)$ of $\CC[M(\G_c)]$ which is a part of a basis of
$\CC[M(\G_c)]$.}
\nl
(See \cite{L19},\cite{L23}).
We have $Con_c^+\sub S(\G_c)$. More precisely, (a)
restricts to a bijection of
$$\bX_{\G_c,*}:=\{(\G',\G'')\in\bX_{\G_c};\G'=\G''\}$$
onto $Con_c^+$. Let
$$\un\bX_{\G_c}=\{\G';(\G',\G')\in\bX_{\G_c,*}\}.$$

In \cite{L20} a bijection $\e$ from a certain basis
of $\CC[M(\G_c)]$ (containing $S(\G_c)$) to $M(\G_c)$ is
defined. This restricts to an injective map from
$S(\G_c)$) to $M(\G_c)$ whose image is equal
to the image of $c\sub M(\G_c)$ (if $c$ is not as in 0.1(a))
and is equal to the image of $c\sub M(\G_c)$, disjoint
union with a single element $(1,?)\in M(\G_c)-c$
(if $c$ is as in 0.1(a)). From the definition of $\e$, the
following holds.

(b) {\it $\e(B)$ appears with nonzero coefficient in $B$
for any $B\in S(\G_c)$.}
\nl
For $\G'_0\in\un\bX_{\G_c}$ let $\bX^{\G'_0}_{\G_c}$
be the set of all $(\G',\G'')\in\bX_{\G_c}$ such that $\G'$
is conjugate to $\G'_0$.
The following statement will be verified in \S1,\S2.

(c) {\it For $\G'_0\in\un\bX_{\G_c}$, the function
$(\G',\G'')\m|\G''|$ on $\bX^{\G'_0}_{\G_c}$ reaches its
maximum at a unique $(\G',\G'')$.}
\nl
For $\G'_0\in\un\bX_{\G_c}$ (with corresponding
$(\G',\G'')$ defined by (c)) we set

$B_{\G'_0}=\ss_{\G',\G''}(1,1)\in S(\G_c)$,

$E_{\G'_0}=$ element of $c$ which maps to $\e(B_{\G'_0})$
under $c\sub M(\G_c)$. (If $c$ is as in 0.1(a), we
necessarily have $\e(B_{\G'_0})\ne(1,?)$.)

We can now define
$${}^1S(\G_c)=\{B_{\G'_0};\G'_0\in\un\bX_{\G_c}\}
\sub S(\G_c),$$
$$\ca_{\G_c}=\{E_{\G'_0};\G'_0\in\un\bX_{\G_c}\}\sub c.$$
The elements of $\ca_{\G_c}$ are said to be the
{\it almost special} representations of $W$ in $c$.

We have a bijection $\ca_{\G_c}@>\si>>Con^+_c$ given by
$E_{\G'_0}\m\ss_{\G'_0,\G'_0}(1,1)$ for
$\G'_0\in\un\bX_{\G_c}$.

We have a bijection ${}^1S(\G_c)@>>>Con^+_c$ given by
$B_{\G'_0}\m\ss_{\G'_0,\G'_0}(1,1)$ for
$\G'_0\in\un\bX_{\G_c}$.

We have a bijection ${}^1S(\G_c)@>>>\ca_{\G_c}$ given by
$B_{\G'_0}\m E_{\G'_0}$ for $\G'_0\in\un\bX_{\G_c}$.

The following statement will be verified in \S1,\S2.

(d) {\it Let $\G'_0\in\un\bX_{\G_c}$. Let
$(\G',\G'')\in\bX^{\G'_0}_{\G_c}$. Assume that
$(x,\s)\in M(\G_c)$ appears with nonzero coefficient
in $\ss_{\G',\G''}(1,1)$. Then $(x,\s)\in M(\G_c)$ appears
with nonzero coefficient in $\ss_{\G'_0,\G'_0}(1,1)$.}
\nl
Assume now that $(x,\s)\in M(\G_c)$
corresponds to $E_{\G'_0}$ under $c\sub M(\G_c)$.
By (b), $(x,\s)$ appears with  nonzero coefficient
in $\ss_{\G',\G''}(1,1)$ where $(\G',\G'')$ is defined by
$\G'_0$ as in (c). Using (d) we see that

(e) {\it $(x,\s)$ appears with  nonzero coefficient
in $\ss_{\G'_0,\G'_0}(1,1)$.}

\subhead 0.3. Notation\endsubhead
Let $F$ be the field $\ZZ/2$. For $a,b$ in $\ZZ$ we write
$a\ll b$ whenever $b-a\ge2$. For $a,b$ in $\ZZ$ let
$[a,b]=\{z\in\ZZ;a\le z\le b\}$. For a finite set $E$
we write $|E|$ for the cardinal of $E$.

\head 1. Classical types\endhead
\subhead 1.1\endsubhead
Let $D\in2\NN$. Let $\ci_D$ be the set of all intervals
$I=[a,b]$ where $1\le a\le b\le D$.
For $I=[a,b],I'=[a',b']$ in $\ci_D$ we write $I\prec I'$
whenever $a'<a\le b<b'$; we write $I\spa I'$
if $a'-b\ge2$ or $a-b'\ge2$. Let $\ci^1_D$ be the set of all
$I=[a,b]\in\ci_D$ such that $a=b\mod2$.
For $I=[a,b]\in\ci^1_D$ we define $\k(I)\in\{0,1\}$ by
$\k(I)=0$ if $a,b$ are even, $\k(I)=1$ if $a,b$ are odd.

A sequence $I_*=(I_1,I_2,\do,I_r)$ in $\ci^1_D$ is said to
be {\it admissible}
if $$I_*=([a_1,b_1],[a_2,b_2],\do,[a_r,b_r]),r\ge1$$
where $$a_1\le b_1,a_2\le b_2,\do,a_r\le b_r,$$
$$a_2-b_1=2,a_3-b_2=2,\do,a_r-b_{r-1}=2.$$
For such $I_*$ we define $\k(I_*)\in\{0,1\}$ by
$\k(I_*)=0$ if all (or some) $a_i,b_i$ are even,
$\k(I_*)=1$ if all (or some) $a_i,b_i$ are odd.

For $I=[a,b]\in\ci^1_D$ let
$$I^{ev}=\{x\in I;x=a+1\mod2\}=\{x\in I;x=b+1\mod2\}.$$

\subhead 1.2\endsubhead
Let $R^1_D$ be the set whose elements are the subsets of $\ci^1_D$.
Let $B\in R^1_D$. We consider the following properties
$(P_0),(P_1)$ that $B$ may or may not have:

$(P_0)$ If $I\in B,I'\in B$ then either $I=I'$ or $I\spa I'$ or $I\prec I'$ or $I'\prec I$.

$(P_1)$ If $I\in B$ and $x\in I^{ev}$
then there exists $I'in B$ such that $x\in I',I'\prec I$.

Let $S_D$ be the set of all $B\in R^1_D$ that satisfy
$(P_0),(P_1)$. (In \cite{L19}, two sets $S_D,S'_D$ are introduced
and showed to be equal. What we call $S_D$ in this paper was
called $S'_D$ in \cite{L19}.)

For $B\in S_D,I\in B$
let $m_{I,B}=|\{I'\in B;I\sub I'\}|$. 

\subhead 1.3\endsubhead
For $B\in S_D$, $I=[a,b]\in B$ we set
$$X_{I,B}=\{I'\in B;I'\prec I,m_{I',B}=m_{I,B}+1\}.$$
Assuming that $a<b$, we show:

(a) {\it $X_{I,B}$ is an admissible sequence
$$([a_1,b_1],[a_2,b_2],\do,[a_r,b_r])$$ (see 1.1) with
$a_1=a+1,b_r=b-1$. Moreover $\k(X_{I,B})=1-\k(I)$.}
\nl
We have $a+1\in I^{ev}$. By $(P_1)$ we can find
$[a_1,b_1]\in B$ such that $a<a_1\le a+1\le b_1<b$; we must have
$a_1=a+1$ and we can assume that $b_1$ is maximum possible.
Then $m_{[a_1,b_1],B}=m_{I,B}+1$. If $b_1=b-1$ then we stop.
Assume now that $b\ge b_1+3$. Let $a_2=b_1+2$. We have
$a_2\in I^{ev}$ hence by $(P_1)$ we can find
$[x,b_2]\in B$ such that $a<x\le a_2\le b_2<b$; we can
assume
that $b_2$ is maximum possible. Then $m_{[x,b_2],B}=m_{I,B}+1$.
Since $[a_1,b_1]\spa[x,b_2]$, we must have $x=a_2$. If
$b_2=b-1$
then we stop. Assume now that $b\ge b_2+3$. Let $a_3=b_2+2$.
We have $a_3\in I^{ev}$ hence by $(P_1)$ we can find
$[x,b_3]\in B$ such that $a<x\le a_3\le b_3<b$; we can
assume
that $b_3$ is maximum possible. Then $m_{[x,b_3],B}=m_{I,B}+1$.
Since $[a_2,b_2]\spa[x,b_3]$, we must have $x=a_3$.
This process continues in this way and it eventually stops.
This proves (a). (The last statement of (a) is obvious.)

\subhead 1.4\endsubhead
Let $B\in S_D$. For $I=[a,b]\in B$ we have

(a) $|\{I'\in B;I'\sub I\}|=(b-a+2)/2$.
\nl
See \cite{L20, 1.3(d)}. We now write the various $I\in B$ such that $m_{I,B}=1$ in a sequence
$[a_1,b_1],[a_2,b_2],\do,[a_r,b_r]$
where
$$1\le a_1\le b_1\ll a_2\le b_2\ll\do\ll a_r\le b_r.$$
From (a) we deduce
$$\align&|B|=\sum_{i\in[1,r]}(b_i-a_i+2)/2\\&
=-a_1/2 -((a_2-b_1)+(a_3-b_2)+\do+(a_r-b_{r-1}))/2
+b_r/2+r\\&\le
(b_r-a_1)/2-(r-1)+r\le(D-1)/2+1=(D+1)/2.\endalign$$
Since $D\in2\NN$ it follows that $|B|\le D/2$. We see that:

(b) {\it the condition that $|B|=D/2$ is that either}

-{\it each of $a_2-b_1,a_3-b_2,\do,a_r-b_{r-1}$ equals $2$
except one of them which equals $3$ and $a_1=1,b_r=D$,} or

-{\it each of $a_2-b_1,a_3-b_2,\do,a_r-b_{r-1}$ equals $2$
and $a_1=1,b_r=D-1$}, or

-{\it each of $a_2-b_1,a_3-b_2,\do,a_r-b_{r-1}$ equals $2$
and $a_1=2,b_r=D$.}
\nl
We set
$S_D^{D/2}=\{B\in S_D;|B|=D/2\}$.

\subhead 1.5\endsubhead
Let ${}^1S_D$ be the set of all $B\in S_D$
such that $\k(I)=1$ for any $I\in B$ with $m_{I,B}=1$.

Let $B\in S_D$. We set
${}^1B=B-\{I\in B;\k(I)=0,m_{I,B}=1\}\in R^1_D$.
From the definitions it is clear that ${}^1B\in S_D$.
We show that

(a) ${}^1B\in{}^1S_D$.
\nl
Let $I'\in{}^1B$ be such that $m_{I',{}^1B}=1$.
We have either $m_{I',B}=1,\k(I')=1$ or else
$m_{I',B}=2$ and $I'\prec I''$ for some $I''\in B$ with $m_{I'',B}=1$, $\k(I'')=0$. In the second case we have
$I'\in X_{I'',B}$, so that from 1.3(a) we have $\k(I')=1$.
This proves (a).

Thus $B\m{}^1B$ is a well defined map $S_D@>>>{}^1S_D$.

\subhead 1.6\endsubhead
Let $B\in{}^1S_D$. We write the various $I\in B$ such that
$m_{I,B}=1$ in a sequence
$$[a^1_1,b^1_1],[a^1_2,b^1_2],\do,[a^1_{r_1},b^1_{r_1}],$$
$$[a^2_1,b^2_1],[a^2_2,b^2_2],\do,[a^2_{r_2},b^2_{r_2}],$$
$$\do,$$
$$[a^s_1,b^s_1],[a^s_2,b^s_2],\do,[a^s_{r_s},b^s_{r_s}]$$
whose first $r_1$ terms form an admissible sequence $I_{*1}$,
the next $r_2$ terms form an admissible sequence $I_{*2}$,
$\do$, and the last $r_s$ terms form an admissible sequence
$I_{*s}$; we also assume that
$$a^2_1\ge b^1_{r_1}+4,a^3_1\ge b^2_{r_2}+4,\do,
a^s_1\ge b^{s-1}_{r_{s-1}}+4.$$
Here we have
$$r_1\ge1,r_2\ge1,\do,r_s\ge1,s\ge0,
\k(I_{*1})=1,\k(I_{*2})=1,\do,\k(I_{*s})=1.$$

Let $Z(B)$ be the subset of $\ci^1_D$ consisting of:

all $[a^i_1-1,b^i_{r_i}+1]$ ($i\in[1,s]$) such that
$a^i_1\ge2$ (this is automatic if $i\ge2$);

all $[u,u]$ with $u$ even, $b^{i-1}_{r_{i-1}}+1<u<a^i_1-1$
for some $i\in[2,s]$ (if $s>1$);

all $[u,u]$ with $u$ even, $1<u<a^1_1-1$ (if $s>0$);

all $[u,u]$ with $u$ even, $b^s_{r_s}+1<u\le D$ (if $s>0$);

all $[u,u]$ with $u$ even, $1<u\le D$ (if $s=0$).

For any subset $U\sub Z(B)$ we set
$B_U=B\sqc U$; then $B_U\in S_D$ and $U\m B_U$ defines a
bijection from the set of subsets of $Z(B)$ to the fibre
at $B$ of the map $S_D@>>>{}^1S_D$, $B'\m{}^1B'$. Note that
$B_\emp=B$ and $B_{Z(B)}\in S_D^{D/2}$. Moreover,
$B\m B_{Z(B)}$ is the bijection ${}^1S_D@>\si>>S_D^{D/2}$
whose inverse is the restriction to $S_D^{D/2}$ of
$S_D@>>>{}^1S_D$, $B'\m{}^1B'$. (We use 1.4(b).)

\subhead 1.7\endsubhead
A subset $B$ of $R^1_D$ is said to be in ${}^1\dot{S}_D$
if it satisfies $(P_0)$ and if each $I\in B$ satisfies
$\k(I)=1$. For $B\in{}^1S_D$ we set
$\dot{B}=\{I\in B;\k(I)=1\}$. Then $B\m\dot{B}$ is a map

(a) ${}^1S_D@>>>{}^1\dot{S}_D$.
\nl
We show:

(b) {\it The map (a) is a bijection.}
\nl
Let $C\in{}^1\dot{S}_D$. For $I\in C$ we set
$\dot{m}_{I,C}=|\{I'\in C;I\sub I'\}|$.

For $k\in\{1,2,3,\do\}$ we set
$C[k]=\{I\in C;\dot{m}_{I,C}=k\}$.

Let $I=[a,b]\in C[k]$. As in 1.6 we can write the
intervals $\{I'\in C[k+1];I'\prec I\}$ in a sequence
$$[a^1_1,b^1_1],[a^1_2,b^1_2],\do,[a^1_{r_1},b^1_{r_1}],$$
$$[a^2_1,b^2_1],[a^2_2,b^2_2],\do,[a^2_{r_2},b^2_{r_2}],$$
$$\do,$$
$$[a^s_1,b^s_1],[a^s_2,b^s_2],\do,[a^s_{r_s},b^s_{r_s}]$$
whose first $r_1$ terms form an admissible sequence $I_{*1}$,
the next $r_2$ terms form an admissible sequence $I_{*2}$,
$\do$, and the last $r_s$ terms form an admissible sequence
$I_{*s}$; we also assume that
$$a^2_1\ge b^1_{r_1}+4,a^3_1\ge b^2_{r_2}+4,\do,
a^s_1\ge b^{s-1}_{r_{s-1}}+4.$$
Here we have
$$r_1\ge1,r_2\ge1,\do,r_s\ge1,s\ge0,
\k(I_{*1})=1,\k(I_{*2})=1,\do,\k(I_{*s})=1.$$
Moreover we have $a^i_1\ge a+2,b^i_{r_i}\le b-2$ for all $i$.

Let $Y_I$ be the subset of $\ci^1_D$ consisting of:

all $[a^i_1-1,b^i_{r_i}+1]$ ($i\in[1,s]$);

all $[u,u]$ with $u$ even, $b^{i-1}_{r_{i-1}}+1<u<a^i_1-1$
for some $i\in[2,s]$ (if $s>1$);

all $[u,u]$ with $u$ even, $a<u<a^1_1-1$ (if $s>0$);

all $[u,u]$ with $u$ even, $b^s_{r_s}+1<u<b$ (if $s>0$);

all $[u,u]$ with $u$ even, $a<u<b$ (if $s=0$).

For $l\ge1$ we set
$B[2l-1]=C[l],B[2l]=\sqc_{I\in C[l]}Y_I$. We set
$B=\sqc_{l\in\{1,2,3,\do\}}B[l]$.
From the definition we see that $B\in{}^1S_D$ and that
$C\m B$ is an inverse of the map ${}^1S_D@>>>{}^1\dot{S}_D$,
$B\m\dot{B}$. This proves (b).

\mpb

We shall view any element $C\in{}^1\dot{S}_D$ as a tableau with
columns indexed by $[1,D]$, with rows indexed by
$\{1,2,3,\do\}$ and with entries in  $\cup_j[a_j,b_j]$.
Any entry in the $s$-column is equal to $s$; the $k$-th
row consists of the elements in $\cup_{I\in C[k]}I$.

\subhead 1.8\endsubhead
Let $C\in{}^1\dot{S}_D$. It is an unordered set of
intervals $[a_1,b_1],[a_2,b_2],\do,[a_t,b_t]$.
We can order them by requiring that $b_1<b_2<\do<b_t$.
We view $C$ as a tableau as in 1.7.
We associate to $C$ a new tableau $\ddot{C}$ with
columns indexed by $[1,D]$, with rows indexed by
$\{1,2,3,\do\}$ and with entries in $\cup_j[a_j,b_j]$.
This is obtained by moving the entry of $C$ in
the $s$-column and row $k$ to the same $s$-column and
to row $k+j$ where $j\in[0,t-1]$ is defined by
$b_j<s\le b_{j+1}$  (with the convention $b_0=0$); note
that $s\le b_t$ whenever the $s$-column of $C$ is nonempty.

For example, ${}^1\dot{S}_4$ consists of $5$ tableaux:
$(\emp)$, $(1)$; $(3)$; $\mat1&{}&3\endmat$; $(123)$. The corresponding
tableaux $\ddot{C}$
are $(\emp)$, $(1)$; $(3)$; $\mat 1&{}&{}\\{}&{}&3\endmat$; $(123)$. 
Now ${}^1\dot{S}_6$ consists of $14$ tableaux:

$(\emp)$; $(1)$; $(3)$; $(5)$; $\mat1&{}&3\endmat$;
$\mat3&{}&5\endmat$;
$\mat1&{}&{}&{}&5\endmat$; $\mat1&{}&3&{}&5\endmat$;
$(123)$; $\mat123&{}&5\endmat$; $(345)$; $\mat1&{}&345\endmat$;
$(12345)$; $\mat 1&2&3&4&5\\{}&{}&3&{}&{}\endmat$.

The corresponding tableaux $\ddot{C}$ are

$(\emp)$; $(1)$; $(3)$; $(5)$;
$\mat1&{}&{}\\{}&{}&3\endmat$;
$\mat3&{}&{}\\{}&{}&5\endmat$;
$\mat1&{}&{}&{}&{}\\{}&{}&{}&{}&5\endmat$;
$\mat1&{}&{}&{}&{}\\{}&{}&3&{}&{}\\{}&{}&{}&{}&5\endmat$;
$(123)$; $\mat123&{}\\{}&5\endmat$; $(345)$;
$\mat1&{}\\{}&345\endmat$;
$(12345)$; $\mat 1&2&3&4&5\\{}&{}&3&{}&{}\endmat$.

Here are some further examples.

If $C=\mat1&2&3&4&5\\ {}&{}&3&{}&{}\endmat$
then $\ddot{C}=\mat 1&2&3&{}&{}\\ {}&{}&3&4&5\endmat$.

\mpb

If $C=\mat 1&2&3&4&5&6&7\\{}&{}&3&{}&{}&{}&{}\endmat$ 
then $\ddot{C}=\mat 1&2&3&{}&{}&{}&{}\\{}&{}&3&4&5&6&7\endmat$.

\mpb

If $C=\mat 1&2&3&4&5&6&7\\{}&{}&{}&{}&5&{}&{}\endmat$ then
$\ddot{C}=\mat 1&2&3&4&5&{}&{}\\{}&{}&{}&{}&5&6&7\endmat$.

\mpb

If $C=\mat1&2&3&4&5&6&7\\{}&{}&3&{}&5&{}&{}\endmat$ then
$\ddot{C}=\mat1&2&3&{}&{}&{}&{}\\{}&{}&3&4&5&{}&{}
\\{}&{}&{}&{}&5&6&7\endmat$.

\mpb

If $C=\mat1&2&3&4&5&6&7\\{}&{}&3&4&5&{}&{}\endmat$ then
$\ddot{C}=\mat1&2&3&4&5&{}&{}\\{}&{}&3&4&5&6&7\endmat$.

\mpb

If $C=\mat1&2&3&4&5&6&7&8&9\\{}&{}&3&4&5&6&7&{}&{}
\\{}&{}&{}&{}&5&{}&{}&{}&{}\endmat$
then

$\ddot{C}=\mat1&2&3&4&5&{}&{}&{}&{}
\\{}&{}&3&4&5&6&7&{}&{}\\{}&{}&{}&{}&5&6&7&8&9\endmat$.

We show: 

(a) {\it Let $j\in[1,t]$. Let $k$ be such that
 $[a_j,b_j]\in C[k]$. In row $j$ of $\ddot{C}$, $b_j$
 appears and $b_j+1$ does not appear.}
 \nl
In rows $j+1,j+2,\do,j+k-1$ of $\ddot{C}$, $b_j$ and $b_j+1$
appear. In any other row of $\ddot{C}$, $b_j$ and $b_j+1$ do not appear.
\nl
 Assume first that $b_j<D$.
 Then in $C$, $b_j$ appears in rows
$1,2,\do,k$ and $b_j+1$ appears in rows $1,2,\do,k-1$.
Since $b_{j-1}<b_j\le b_j$, $b_j<b_j+1\le b_{j+1}$ we see
that in $\ddot{C}$, $b_j$ appears in rows
$1+(j-1),2+(j-1),\do,k+(j-1)$ and $b_j+1$ appears in rows
$1+j,2+j,\do,(k-1)+j$. This proves (a) in our case. Now
assume that $b_j=D$ (in this case $j=t$ and $k=1$). Then
in $C$, $b_t$ appears in row $1$ and in no other row. We
have $b_{t-1}<b_t\le b_t$. Hence in $\ddot{C}$, $b_t$ appears in
row $1+(t-1)$ and in no other row. Thus (a) again holds.

We show: 

(b) {\it Let $i\in[1,t]$. Let $k$ be such that
 $[a_i,b_i]\in C[k]$. Define $j\in[0,t-1]$ by
$b_j<a_i\le b_{j+1}$. In row $k+j$ of $\ddot{C}$, $a_i$ appears
and $a_i-1$ does not appear.}
\nl
In rows $j+1,j+2,\do,j+k-1$ of $\ddot{C}$, $a_i$ and $a_i-1$
appear. In any other row of $\ddot{C}$, $a_i$ and $a_i-1$ do not
appear.
\nl
Assume first that $a_i>1$. Then in $C$, $a_i$ appears in rows
$1,2,\do,k$ and $a_i-1$ appears in rows $1,2,\do,k-1$.
Then
(since $b_j,a_i$ are odd) we have $b_j<a_i-1\le b_{j+1}$.
Hence in $\ddot{C}$, $a_i$ appears in rows $1+j,2+j,\do,k+j$
and $a_i-1$ appears in rows $1+j,2+j,\do,(k-1)+j$.
This proves (b) in our case. Now assume that $a_i=1$
(in this case we have $k=1$). Then in $C$, $a_i$ appears
in row $1$ and in no other row. We have $b_0<a_i\le b_1$.
Hence in $\ddot{C}$, $a_i$ appears in row $1$ and in no other
row. Thus (b) again holds.

Now let $h\in\cup_j[a_j,b_j]$ be such that
$h\ne a_j,h\ne b_j$ for all $j\in[1,t]$. We show:

(c) {\it Any row of $\ddot{C}$ that contains $h$ must
also contain $h+1$.}

(d) {\it Any row of $\ddot{C}$ that contains $h$ must
also contain $h-1$.}
\nl
There is a well defined $j\in[0,t-1]$ such that
$b_j<h<b_{j+1}$.

We prove (c). Assume first that $h+1<b_{j+1}$.
Then in $C$, $h$ appears in rows
$1,2,\do,k$ and $h+1$ appears in rows $1,2,\do,k$
(for some $k$).
In $\ddot{C}$, $h$ appears in rows $j+1,j+2,\do,j+k$ and
$h+1$ appears in rows $j+1,j+2,\do,j+k$. Hence
in this case (c) holds.
Next we assume that $h+1=b_{j+1}$.
Then in $C$, $h$ appears in rows
$1,2,\do,k$ and $h+1$ appears in rows
in rows $1,2,\do,k+1$ (for some $k$).
In $\ddot{C}$, $h$ appears in rows $j+1,j+2,\do,j+k$ and
$h+1$ appears in rows $j+1,j+2,\do,j+k+1$.
We see again that (c) holds.

We prove (d). Assume first that $b_j<h-1$. Then in $C$, $h$
appears in rows
$1,2,\do,k$ and $h-1$ appears in rows $1,2,\do,k$
(for some $k$).
In $\ddot{C}$, $h$ appears in rows $j+1,j+2,\do,j+k$ and
$h-1$ appears in rows $j+1,j+2,\do,j+k$. Hence
in this case (d) holds.

Next we assume that $b_j=h-1$. Then in $C$, $h$ appears in
rows $1,2,\do,k$ and $h-1$ appears in rows $1,2,\do,k,k+1$
(for some $k$). Moreover in $\ddot{C}$, $h$ appears in rows
$j+1,j+2,\do,j+k$ and $h-1$ appears in rows
$j,j+1,j+2,\do,j+k$. We see again that (d) holds.

From (a)-(d) we deduce:

(e) {\it For $j\in[1,t]$, the row $j$ of $\ddot{C}$
consists of
$a_{i_j},a_{i_j}+1,a_{i_j}+2,\do,b_j$ for a well defined
$i_j\in[1,t]$ such that $a_{i_j}\le b_j$. Moreover,
$j\m i_j$ is a permutation of $[1,t]$.}

We show:

(f) {\it For $u\in[2,t]$ we have $a_{i_{u-1}}<a_{i_u}$.}
\nl
We set $i=i_u$.

Assume first that $[a_i,b_i]\in C[k],k\ge2.$
By (b), one row of $\ddot{C}$ contains $a_i$ but not $a_i-1$
(hence it is necessarily the row $u$) and the row just
above it (that is row $u-1$) contains $a_i$ and $a_i-1$.
(We use that $k\ge2$.) Now that row consists of
$a_{i_{u-1}},a_{i_{u-1}}+1,\do,b_{u-1}-1,b_{u-1}$.
Thus we have $ a_{i_{u-1}}\le a_i-1<a_i\le b_{u-1}$.
In particular, $a_{i_{u-1}}<a_i$.

Next we assume that $[a_i,b_i]\in C[1]$.
Now $[a_i,b_u]$ is contained in the union of all $I\in C[1]$
and consists of consecutive numbers. Hence it is contained
in one such $I$ which is necessarily $[a_i,b_i]$. Thus we
have $[a_i,b_u]\sub[a_i,b_i]$.

Assume that $a_i\le b_{u-1}$.
In row $1$ of $C$ we have $a_i\le b_{u-1}<b_u$.
In row $2$ of $C$, $a_i$ is missing. Since $a_i\le b_{u-1}$
the unique entry $a_i$ in $\ddot{C}$ appears in a row $\le u-1$.
In particular the row $u$ of $\ddot{C}$ does not contain $a_i$;
but it contains $b_u$. This contradicts $a_i=a_{i_u}$. We
see that we must have $b_{u-1}<a_i$. But
$a_{i_{u-1}}\le b_{u-1}$ hence $a_{i_{u-1}}<a_i$. This
proves (f).

\subhead 1.9\endsubhead
Let ${}^1\ddot{S}_D$ be the set of tableaux with
columns indexed by $[1,D]$, with rows indexed by
$\{1,2,3,\do\}$ and with entries in  $[1,D]$
such that
any entry in the $s$-column is equal to $s$; for any
$k\in[1,t]$ (some $t$), the row $k$ consists of the elements in some $I_k=[c_k,d_k]\in\ci^1_D$ with $\k(I_k)=1$; for
$k>t$ the row $k$ contains no entries; we have
$c_1<c_2<\do<c_t,d_1<d_2<\do<d_t$.

For $X,c_1<c_2<\do<c_t,d_1<d_2<\do<d_t$ as above
we define a tableau $\dot{X}$
with columns indexed by $[1,D]$, with rows indexed by
$\{1,2,3,\do\}$ and with entries in $\cup_j[c_j,d_j]$. This
is obtained by moving the entry of $X$ in the $s$-column
and row $k$ to the same $s$-column and to row $k-j$ where
$j\in[0,t-1]$ is defined by $d_j<s\le d_{j+1}$  (with the
convention $d_0=0$); note that we necessarily have $k>j$.
(Indeed, we have $s\le d_k$; if $k\le j$ then $d_k\le d_j$,
hence $s\le d_j$, contradicting $d_j<s$.)

From the definitions we see that $\dot{X}\in{}^1\dot{S}_D$ and that
$X\m\dot{X}$ is a bijection ${}^1\ddot{S}_D@>>>{}^1\dot{S}_D$ inverse to
$C\m\ddot{C}$, ${}^1\dot{S}_D@>>>{}^1\ddot{S}_D$.

\subhead 1.10\endsubhead
Let $U_D$ be the set of all tableaux
$$\mat c_1&c_2&\do&c_t\\d'_1&d'_2&\do&d'_t\endmat$$
where $c_1<c_2<\do<c_t$ are odd integers in $[1,D]$,
$d'_1<d'_2<\do<d'_t$ are even integers in $[1,D]$
and $c_1<d'_1,c_2<d'_2,\do,c_t<d'_t$.

We have an obvious bijection ${}^1\ddot{S}_D@>\si>>U_D$,
$$\align&(X,c_1<c_2<\do<c_t,d_1<d_2<\do<d_t)\m\\&
\mat c_1&c_2&\do&c_t\\d_1+1&d_2+1&\do&d_t+1\endmat.
\endalign$$

\subhead 1.11\endsubhead
Let $\ddot{\Si}_D$ be the set of all symbols
$$\L=\mat i_1&i_2&\do&i_{(D+2)/2}\\ j_1&j_2&\do&j_{(D+2)/2}
\endmat$$
where
$$\{i_1,i_2,\do,i_{(D+2)/2}\}\sqc
\{j_1,j_2,\do,j_{(D+2)/2}\}=[0,D+1],$$
$$i_1<i_2<\do<i_{(D+2)/2},j_1<j_2<\do<j_{(D+2)/2},$$
$$i_1<j_1,i_2<j_2,\do,i_{(D+2)/2}<j_{(D+2)/2}.$$

We then have $i_1=0,j_{(D+2)/2}=D+1$.

For $\L$ as above let $c_1<c_2<\do<c_t$ be the odd numbers
in $\{i_1,i_2,\do,i_{(D+2)/2}\}$ (in increasing order)
and let $d'_1<d'_2<\do<d'_{t'}$ be the even numbers
in $\{j_1,j_2,\do,j_{(D+2)/2}\}$ (in increasing order).
We have necessarily $t=t'$. We show:
$$c_1<d'_1,c_2<d'_2,\do,c_t<d'_t.\tag a$$

Assume now that for some $s\in[0,t],s<t$ we already know
that $c_1<d'_1,c_2<d'_2,\do,c_s<d'_s$. We show that
$c_{s+1}<d'_{s+1}$.

Assume that $d'_{s+1}\le c_{s+1}$. 
Let $Z=\{i_k;k\in[1,(D+2)/2];i_k<d'_{s+1}\}$.
Then
$$Z=\{0,2,4,\do,d'_{s+1}-2\}\sqc\{c_1,c_2,\do,c_s\}-
\{d'_1,d'_2,\do,d'_s\}.$$
(We use that $c_1<d'_1,c_2<d'_2,\do,c_s<d'_s$. We also use
that $d'_{s+1}\le c_{s+1}$.)
We have $|Z|=|(0,2,4,\do,d'_{s+1}-2)|$.
We have $d'_{s+1}=j_m$ for some $m\in[1,(D+2)/2]$ and
$i_m<d'_{s+1}$ that is $i_m\in Z$. It follows that
$\{i_1,i_2,\do,i_m\}\sub Z$ so that $m\le|Z|$.
Let $$Z'=\{j_k;k\in[1,(D+2)/2];j_k\le d'_{s+1}\}.$$
We have $|Z'|=m$.
Now
$$Z'=\{1,3,5,\do,d'_{s+1}-1\sqc
\{d'_1,d'_2,\do,d'_s,d'_{s+1}\}-\{c_1,s_2,\do,c_s\}$$
so that $|Z'|=|(1,3,5,\do,d'_{s+1}-1)+1)|$. Since
$|Z'|=m$ and $m\le Z$ we have $|Z'|\le|Z|$ so that
$$|(1,3,5,\do,d'_{s+1}-1)+1)|\le|(0,2,4,\do,d'_{s+1}-2)|.$$
This is obviously not true. This proves (a).

\mpb

From (a) we see that
$$\L\m\mat c_1&c_2&\do&c_t\\d'_1&d'_2&\do&d'_t\endmat$$
(as described above) defines a map
$$\ddot{\Si}_D@>>>U_D.\tag b$$
We show:

(c) {\it The map (b) is injective.}
\nl
To any
$$\mu=\mat c_1&c_2&\do&c_t\\d'_1&d'_2&\do&d'_t\endmat\in U_D
$$
we associate a sequence
$$\mu'=(i_1,i_2,\do,i_{(D+2)/2})$$
and a sequence
$$\mu''=(j_1,j_2,\do,j_{(D+2)/2})$$
as follows.

$\mu'$ consists of the elements in $\{c_1,c_2,\do,c_t\}$
and those in $\{0,2,4,\do,D\}-\{d'_1,d'_2,\do,d'_t\}$
(in increasing order).

$\mu''$ consists of the elements in $\{d'_1,d'_2,\do,d'_t\}$
and those in $\{1,3,5,\do,D+1\}-\{c_1,c_2,\do,c_t\}$
(in increasing order).

From the definition we see that if
$\L\in\ddot{\Si}_D$ has image $\mu\in U_D$ under (b)
then $\L=\mat \mu'\\ \mu''\endmat$.
From this it is clear that the map (b) is injective.
This proves (c).

\subhead 1.12\endsubhead
We show:

(a) {\it The injective map $\ddot{\Si}_D@>>>U_D$ in 1.11(b) is a
bijection.}
\nl
Note that $\ddot{\Si}_D$ can be viewed as the set of standard Young
tableaux
attached to a partition with two equal parts of size
$(D+2)/2$. The number of such standard tableaux can be
computed from the hook length formula so that it is equal to
$(D+2)!/((D+2)/2)!(D+4)/2)!)$ that is to the Catalan number
$Cat_{(D+2)/2}$. (This interpretation of Catalan numbers in
terms of standard Young tableaux has been known before.)

From the bijections
$$U_D@<<<{}^1\ddot{S}_D@>>>{}^1\dot{S}_D@<<<
{}^1S_D@>>>S_D^{D/2}$$
(see 1.10,1.9,1.7,1.6) we see that $|U_D|=|S_D^{D/2}|$.
By \cite{LS}, $|S_D^{D/2}|$ is equal to the Catalan number
$Cat_{(D+2)/2}$.
We see that the map in (a) satisfies
$|\ddot{\Si}_D|=|U_D|=Cat_{(D+2)/2}$. It follows that this map
is a bijection.

(It is likely that (a) has a more direct proof which does
not rely on \cite{LS}.)

We show:

(b) {\it If $\mu\in U_D$ and if $\mu',\mu''$ are as in
the proof of 1.11(c), then
$\mat \mu'\\ \mu''\endmat\in\ddot{\Si}_D$.
Moreover $\mu\m\mat \mu'\\ \mu''\endmat$ is the inverse of
the map $\ddot{\Si}_D@>>>U_D$ in 1.11.}
\nl
If $\mu\in U_D$, then by (a) we can find $\L\in\ddot{\Si}_D$
whose
image under the map 1.11(b) is $\mu$. By the proof of
1.11(c) we have $\L=\mat \mu'\\ \mu''\endmat$.
It follows that  $\mat \mu'\\ \mu''\endmat\in\ddot{\Si}_D$.

\subhead 1.13\endsubhead
Let $V_D$ be the $F$-vector space with basis $e_1,e_2,\do,e_D$ and with the symplectic form $(,):V\T V@>>>F$ given by $(e_i,e_j)=1$ if $i-j=\pm1$, $(e_i,e_j)=0$, otherwise. For
any subset $J$ of $[1,D]$ we set $e_J=\sum_{j\in J}e_j\in V_D$.

For $B\in S_D$ let $<B>$ be the subspace of $V_D$
spanned by $\{e_I;I\in B\}$. (This is actually a basis of
$<B>$, see \cite{L19}.)

For $j\in[1,D]$ and $B\in S_D$ we set 
$B_j=\{I\in B;j\in I\}$ and
$$\e_j(B)=|B_j|(|B_j|+1)/2\in F$$.
For $B\in S_D$ we set
$$\e(B)=\sum_{j\in[1,D]}\e_j(B)e_j\in V_D.$$
We show:
$$\e(B)=\sum_{I\in B;m_{I,B}\in2\NN+1}e_I.\tag a$$
An equivalent statement is:

{\it If $j\in[1,D]$ then $|\{I\in B_j,m_{I,B}\in2\NN+1\}|$ is
even if $|B_j|\in(4\ZZ)\cup(4\ZZ+3)$ and is odd if
$|B_j|\in(4\ZZ+1)\cup(4\ZZ+2)$.}
\nl
This follows immediately from the following statement
(which holds by the definition of $S_D$):

{\it $B_j$ consists of intervals
$I_k\prec I_{k-1}\prec\do\prec I_1$ in $\ci^1_D$ such that
$m_{I_k,B}=k,m_{I_{k-1},B}=k-1,\do,m_{I_1,B}=1$.}
\nl
For $C\in{}^1\dot{S}_D$ we define
$$\dot{\e}(C)=\sum_{I\in C}e_I.\tag b$$
For $\ddot{C}\in{}^1\ddot{S}_D$ we define
$$\ddot{\e}(\ddot{C})=\sum_ke_{[c_k,d_k]}\in V_D\tag c$$
where $c_k,c_k+1,c_k+2,\do,d_k$ are the entries in the
$k$-th row of $\ddot{C}$.
From the definitions we have
$$\dot{\e}(C)=\ddot{\e}(\ddot{C})\tag d$$
if $C,\ddot{C}$ correspond to each other under the bijection in
1.9.

\subhead 1.14\endsubhead
By \cite{L19, 1.16}, $B\m\e(B)$ is an injective map
$\e:S_D@>\si>>V_D$. By 1.13(a), for $B\in{}^1S_D$ we have
$\e(B)=\dot{\e}(\dot{B})$ ($\dot{\e}$ as in 1.13(b)). Hence the
restriction of $\e$ to ${}^1S_D$ can be identified with
$\dot{\e}:{}^1\dot{S}_D@>>>V_D$ via the bijection 1.7(b).
In particular, $\dot{\e}$ is injective.
Using 1.13(d) we see that via the bijection in 1.9,
$\dot{\e}:{}^1\dot{S}_D@>>>V_D$ becomes $\ddot{\e}:
{}^1\ddot{S}_D@>>>V_D$. In
particular, $\ddot{\e}$ is injective.

\subhead 1.15\endsubhead
Let $\Si_D$ be the set of all unordered pairs
$\mat A\\B\endmat$ of subsets of $[0,D+1]$ such that
$[0,D+1]=A\sqc B$, $|A|=|B|\mod4$.

There is a unique bijection $f:V_D@>>>\Si_D$ such that
$$f(0)=\mat 0&2&4&\do&D\\1&3&5&\do& D+1\endmat$$
and such that
if $x\in V_D$, $f(x)=\mat A\\ B\endmat$ and $i\in[1,D]$ then
$$f(x+e_i)=\mat A\sharp\{i,i+1\}\\ B\sharp\{i,i+1\}\endmat$$
where $\sharp$ is the symmetric difference; it
follows that for $1\le i<j\le D$ we have
$$f(x+e_i+e_{i+1}+\do+e_j)=
\mat A\sharp\{i,j+1\}\\B\sharp\{i,j+1\}\endmat.$$

\subhead 1.16\endsubhead
We can regard $\ddot{\Si}_D$ as a subset of $\Si_D$.
If $\ddot{C}\in{}^1\ddot{S}_D$ corresponds to $\mu\in U_D$
under 1.10 then from the definitions we have
$f(\ddot{\e}(\ddot{C}))=\mat \mu'\\ \mu''\endmat$ (notation of
1.12(b)). In particular we have

(a) $f(\ddot{\e}(\ddot{C}))\in\ddot{\Si}_D$
\nl
and (using 1.12(b)) we see that

(b) $\ddot{C}\m f(\ddot{\e}(\ddot{C}))$
{\it is a bijection}
${}^1\ddot{S}_D@>\si>>\ddot{\Si}_D$.

\subhead 1.17\endsubhead
We have $V_D=V_D^0\op V_D^1$ where $V_D^0$ is the subspace
spanned by $\{e_2,e_4,\do,e_D\}$ and $V_D^1$ is the
subspace spanned by $\{e_1,e_3,\do,e_{D-1}\}$.
For $I\in\ci^1_D$ we have $I=I^0\sqc I^1$ where
$I^0=I\cap\{2,4,\do,D\}$, $I^1=I\cap\{1,3,\do,D-1\}$.
As shown in \cite{L19}, for $B\in S_D$ we have
$<B>=<B>_0\op<B>_1$ where
$$<B>_0=<B>\cap V^0_D,<B>_1=<B>\cap V^1_D;$$
moreover,

(a) {\it $<B>_1$ has basis $\{e_{I^1};I\in B,\k(I)=1\}$,}

{\it $<B>_0$ has basis $\{e_{I^0};I\in B,\k(I)=0\}$. }

\subhead 1.18\endsubhead
If $\cl$ is a subspace of $V^\d_D$ ($\d\in\{0,1\}$) we set
$$\cl^!=\{x\in V^{1-\d};(x,\cl)=0\}.$$
Let $\cc(V^\d_D)$ be the set of subspaces of $V^\d_D$ of the
form $<B>_\d$ for some $B\in S_D$.
If $\cl\in\cc(V^\d_D)$ we have $\cl^!\in\cc(V^{1-\d}_D)$;
see \cite{L19, \S2}.
Let $\ca(V^1_D)$ be the set of all   
$(\cl,\cl')\in\cc(V^1_D)\T\cc(V^1_D)$
such that $\cl\sub\cl'$ and $\cl\op\cl'{}^!=<B>$
for some $B\in S_D$.

(a) {\it If $B\in S_D$ then $B\m(<B>_1,<B>_0^!)$ is a
bijection $\Ph:S_D@>>>\ca(V^1_D)$;}
\nl
see \cite{L19, \S2}.
Let $\ca_*(V^1_D)$ be the set of all
$(\cl,\cl')\in\cc(V^1_D)\T\cc(V^1_D)$
such that $\cl=\cl'$. In \cite{L19,\S2} it is shown that
$\ca_*(V^1_D)\sub\ca(V^1_D)$; more precisely if
$\cl\in\cc(V^1_D)$ then $\cl\op\cl^!=<B>$ for a well
defined $B\in S_D^{D/2}$. Moreover
$B\m(<B>_1,<B>_0^!)$ is a bijection $S_D^{D/2}@>>>\ca_*(V^1_D)$ and
$(\cl,\cl')\m\cl=\cl'$ is a bijection
$\ca_*(V^1_D)@>>>\cc(V^1_D)$. The composition of these bijections
is a bijection $B\m<B>_1$,

(b) $S_D^{D/2}@>\si>>\cc(V^1_D)$.
\nl
Next we note that $B\m<B>_1$ is also a bijection

(c) ${}^1S_D@>\si>>\cc(V^1_D)$.
\nl
This follows from (b) since the bijection (b) is a composition
$$S_D^{D/2}@>\si>>{}^1S_D@>>>\cc(V^1_D)$$
where the fist map is the bijection $B\m{}^1B$ and the second map is
the map in (c).

Here we use that

(d) $<B>_1=<{}^1B>_1$ for any $B\in S_D$,
\nl
which follows from definitions.

We show:

(e) {\it For any $\cl\in\cc(V^1_D)$, the set
$$\{\cl'\in\cc(V^1_D);(\cl,\cl')\in\ca(V^1_D)\}$$
contains a unique $\cl'$ with $|\cl'|$ maximal.}
\nl
An equivalent statement is:

(f) {\it For any $\cl\in\cc(V^1_D)$ the set
$\{B'\in S_D;<B'>_1=\cl\}$ contains a unique $B'$ with
$|<B'>_0^!|$ maximal (that is $\dim(<B'>_0)|$ minimal).}
\nl
By (b) we have $\cl=<B>_1$ for a well defined $B\in S_D^{D/2}$.
The condition that $<B'>_1=<B>_1$ is equivalent to
$<{}^1B'>_1=<{}^1B>_1$ (see (d)) and this is equivalent to
${}^1B'={}^1B$ (see (c)). Hence the set in (f) is equal to

$\{B'\in S_D;{}^1B'={}^1B\}$.
\nl
By the results in 1.6 this is the same as
$\{({}^1B)_U;U\in Z({}^1B)\}$.
By  1.17(a), for $U\in Z({}^1B)$ we have
$$\dim(({}^1B)_U)_0=\dim({}^1B)_0+|U|.$$
This is $\ge\dim ({}^1B)_0$ with equality if and only if
$U=\emp$. This proves (f) and hence (e).

\mpb

For $\cl\in\cc(V^1_D)$ we denote by $\cl^{max}$ the
element $\cl'$ in (e) with $|\cl'|$ maximal.
Let $\ca^*(V^1_D)$ be the set of all
$(\cl,\cl')\in\cc(V^1_D)\T\cc(V^1_D)$
such that $\cl'=\cl^{max}$. We have 
$\ca^*(V^1_D)\sub\ca(V^1_D)$ and

(g) {\it $(\cl,\cl')\m\cl$ is a bijection
$\ca^*(V^1_D)@>>>\cc(V^1_D)$. }
\nl
From (c),(g) we see (using the proof of (f)) that

(h) {\it $B\m(<B>_1,<B>_0^!)$ is a bijection
${}^1S_D@>\si>>\ca^*(V^1_D)$. }

\subhead 1.19\endsubhead
In this subsection we assume that $D\in\{2,4,6\}$.
In each case we give a table with rows of the form
$\a....\b....\g$
where $\a\in S_D^{D/2}, \b\in {}^1S_D$ corresponds to $\a$
and $\g\in\ddot{\Si}_D$. We write an interval $[a,b]$
as $a,a+1,a+2,\do,b$ (without commas).

\mpb

$D=2$

$\{2\}......\{\emp\}........\mat 0&2\\1&3\endmat$

$\{1\}....\{1\}........ \mat 0&1\\2&3\endmat$

\mpb

$D=4$

\mpb

$\{2,4\}......\{\emp\}......  \mat 0&2&4\\1&3&5\endmat$

$\{1,4\}.....\{1\}.......      \mat 0&1&4\\2&3&5\endmat$

$\{3,234\}...\{3\}......   \mat 0&2&3\\1&4&5\endmat$

$\{1,3\}....\{1,3\}...... \mat 0&1&3\\2&4&5\endmat$

$\{2,123\}....\{2,123\}......\mat 0&1&2\\3&4&5\endmat$

\mpb

$D=6$

\mpb

$\{2,4,6\}......\{\emp\}..... \mat 0&2&4&6\\1&3&5&7\endmat$

$\{1,4,6\}.....\{1\}......    \mat 0&1&4&6\\2&3&5&7\endmat$

$\{3,234,6\}...\{3\}..... \mat 0&2&3&6\\1&4&5&7\endmat$

$\{2,5,456\}...\{5\}.....  \mat 0&2&4&5\\1&3&6&7\endmat$

$\{1,3,6\}....\{1,3\}..... \mat 0&1&3&6\\2&4&5&7\endmat$

$\{1,5,456\}...\{1,5\}......\mat 0&1&4&5\\2&3&6&7\endmat$

$\{3,5,23456\}...\{3,5\}...... \mat0&2&3&5\\1&4&6&7\endmat$

$\{1,3,5\}....\{1,3,5\}....\mat 0&1&3&5\\2&4&6&7\endmat$

$\{2,123,6\}....\{2,123\}...\mat0&1&2&6\\3&4&5&7\endmat$

$\{4,345,23456\}....\{4,345\}.....\mat 0&2&3&4\\1&5&6&7\endmat$

$\{2,4,12345\}....\{2,4,12345\}....
\mat0&1&2&4\\3&5&6&7\endmat$

$\{1,4,345\}....\{1,4,345\}......\mat 0&1&3&4\\2&5&6&7\endmat$

$\{2,123,5\}....\{2,123,5\}......
\mat0&1&2&5\\3&4&6&7\endmat$

$\{3,234,12345\}....\{3,234,12345\}......
\mat0&1&2&3\\4&5&6&7\endmat$

\head 2. Exceptional types\endhead
\subhead 2.1\endsubhead
Let $W,c,\G_c$ be as in 0.1. We must show that 0.2(c),0.2(d)
hold.

If $W$ is of type $A_n,n\ge1$ we have $|c|=1,\G_c=S_1$.
In this case, 0.2(c), 0.2(d) are trivial.

If $W$ is of type $B_n$ or $C_n$, $n\ge2$ or $D_n,n\ge4$,
we can identify $\G_c=V_D^1$ for some $D\in2\NN$.
We can identify $M(\G_c)=V_D$ as in \cite{L19, 2.8(i)}.
In these cases, 0.2(c), 0.2(d) follow from 1.18(e) and the
proof of 1.18(f). Now $\ca_c$ is the same as $\ddot{\Si}_D$ (see
1.11) in the symbol notation \cite{L84} for objects of $\hW$
(assuming that $W$ is of type $D$ and $c$ is a cuspidal family).

\subhead 2.2\endsubhead
In the remainder of this section we assume that $W$ is of
exceptional type. Then we are in one of the following cases.

$|c|=1,\G_c=S_1$;

$|c|=2,\G_c=S'_2$;

$|c|=3,\G_c=S_2$;

$|c|=4,\G_c=S'_3$;

$|c|=5,\G_c=S_3$;

$|c|=11,\G_c=S_4$;

$|c|=17,\G_c=S_5$.

Here $S'_2$ (resp. $S'_3$) is another copy of $S_2$ (resp.
$S_3$).

For the elements $(\G',\G'')\in\bX_{\G_c}$ we use the
notation of \cite{L23}.
Following \cite{L23} we give for each $\G'_0\in\un\bX_{\G_c}$
the list

$L(\G'_0)=\{\G'';(\G',\G'')\in\bX^{\G'_0}_{\G_c}\}$.
\nl
Assume that $|c|=1$. Then $L(S_1)=\{S_1\}$.

Assume that $|c|\in\{2,3\}$.  Then

$L(S_1)=\{S_2,S_1\}$, $L(S_2)=\{S_2\}$.

Assume that $|c|\in\{4,5\}$.  Then

$L(S_1)=\{S_3,S_2,S_1\}$, $L(S_2)=\{S_2\}$,
$L(S_3)=\{S_3\}$.

Assume that $|c|=11$.  Then

$L(S_1)=\{S_4,S_3,S_2S_2,S_2,S_1\}$,
$L(S_2)=\{S_2S_2,S_2\}$,

$L(S_2S_2)=\{\D_8,S_2S_2\}$, $L(S_3)=\{S_3\}$,
$L(\D_8)=\{\D_8\}$, $L(S_4)=\{S_4\}$.

Assume that $|c|=17$.  Then

$L(S_1)=\{S_5,S_4,S_3S_2,S_3,S_2S_2,S_2,S_1\}$,
$L(S_2)=\{S_3S_2,S_2S_2,S_2\}$,

$L(S_2S_2)=\{\D_8,S_2S_2\}$,
$L(S_3)=\{S_3S_2,S_3\}$,
$L(\D_8)=\{\D_8\}$, $L(S_3S_2)=\{S_3S_2\}$,   

$L(S_4)=\{S_4\}$, $L(S_5)=\{S_5\}$.

In each case we see that $L(\G'_0)$ contains a unique term
with $||$ maximum. (It is the first term of $L(\G'_0)$.)
We see that 0.2(c) holds in each case. Now 0.2(d) can be
easily verified using the tables in \cite{L20,\S3}.

\subhead 2.3\endsubhead
Applying $\e$ to $\ss(\G',\G'')$ for each $\G''$ in the
list $L(\G'_0)$ (recall that
$(\G',\G'')\in\bX^{\G'_0}_{\G_c}$) we obtain a list
$L'(\G'_0)$ of elements in $M(\G_c)$; we write in the same
order as the elements of $L(\G'_0)$.
(The notation for elements in $M(\G_c)$ is taken from
\cite{L84}.)

Assume that $|c|=1$. Then $L'(S_1)=\{(1,1)\}$.

Assume that $|c|\in\{2,3\}$.  
Then $L'(S_1)=\{(1,1),(1,\e)\}$, $L'(S_2)=\{(g_2,1)\}$.

Assume that $|c|\in\{4,5\}$.  Then

$L'(S_1)=\{(1,1),(1,r),(1,\e)\}$, $L'(S_2)=\{(g_2,1)\}$,
$L'(S_3)=\{(g_3,1)\}$.

Assume that $|c|=11$.  Then

$L'(S_1)=\{(1,1),(1,\l^1),(1,\s),(1,\l^2),(1,\l^3)\}$,
$L'(S_2)=\{(g_2,1),(g_2,\e'')\}$,

$L'(S_2S_2)=\{(g'_2,1),(g'_2,\e'')\}$,
$L'(S_3)=\{(g_3,1)\}$, $L'(\D_8)=\{(g'_2,\e')\}$,
$L'(S_4)=\{(g_4,1)\}$.

Assume that $|c|=17$.  Then

$L'(S_1)=\{(1,1),(1,\l^1),(1,\nu),(1,\l^2),(1,\nu'),
(1,\l^3),(1,\l^4)\}$,

$L'(S_2)=\{(g_2,1),(g_2,r),(g_2,\e)\}$,
$L'(S_2S_2)=\{(g'_2,1),(g'_2,\e'')\}$,

$L'(S_3)=\{(g_3,1),(g_3,\e)\}$,
$L'(\D_8)=\{(g'_2,\e')\}$, $L'(S_3S_2)=\{(g_6,1)\}$,

$L'(S_4)=\{ (g_4,1)\}$, $L'(S_5)=\{(g_5,1)\}$.

The almost special representations in $c$ are represented by
the first term in each list. They are as follows.

If $|c|=1$ we have $\ca_{\G_c}=\{(1,1)\}$.

If $|c|\in\{2,3\}$ we have $\ca_{\G_c}=\{(1,1),(g_2,1)\}$.

If $|c|\in\{4,5\}$ we have
$\ca_{\G_c}=\{(1,1),(g_2,1),(g_3,1)\}$.

If $|c|=11$ we have
$\ca_{\G_c}=\{(1,1),(g_2,1),(g'_2,1),(g_3,1),
(g'_2,\e'),(g_4,1)\}$.

If $|c|=17$ we have

$\ca_{\G_c}=\{(1,1),(g_2,1),(g'_2,1),(g_3,1),
(g'_2,\e'),(g_6,1),(g_4,1),(g_5,1)\}$.

\subhead 2.4\endsubhead
In the case where $|c|=17$ we have that
$W$ must of type $E_8$.
An element of each list $L'(\G'_0)$ can be identified with
an element of $c$ (under the imbedding $c\sub M(\G_c)$)
represented by its dimension (with
the single exception of $(1,\l^4))$. Then the lists
$L'(\G'_0)$ become:

$L'(S_1)=\{4480,5670,4536,1680,1400,70,?\}$,

$L'(S_2)=\{7168,5600,448\}$,

$L'(S_2S_2)=\{4200,2688\}$,

$L'(S_3)=\{3150,1134\}$,

$L'(\D_8)=\{168\}$, $L'(S_3S_2)=\{2016\}$,

$L'(S_4)=\{1344\}$, $L'(S_5)=\{420\}$.
\nl
Note that the first representation in a given list
$L'(\G'_0)$ has the $b$-invariant (see \cite{L84, (4,1,2)})
strictly smaller than the $b$-invariant of any subsequent
representation in the list. (We expect that this property
holds for any $c$.) This property is similar to the
defining property of special representations \cite{L79a} and
justifies the name of ``almost special'' representations.


\widestnumber\key{ABCD}
\Refs
\ref\key{L79}\by G.Lusztig\paper
 Unipotent representations of a finite Chevalley group of type $E_8$\lb
 \jour Quart.J.Math.\vol30\yr1979\pages315-338\endref
\ref\key{L79a}\by G.Lusztig\paper A class of irreducible representations of a Weyl group\jour Proc.Kon.Nederl. Akad.(A)
\vol82\yr1979\pages323-335\endref
\ref\key{KL}\by D.Kazhdan and G.Lusztig\paper
Representations of Coxeter groups and Hecke algebras\jour
Inv.Math.\vol53\yr1979\pages165-184\endref
\ref\key{L82}\by G.Lusztig\paper A class of irreducible
representations of a Weyl group, II\jour Proc.Kon.Nederl.
Akad.(A)\vol85\yr1982\pages219-226\endref
\ref\key{L84}\by G.Lusztig\book Characters of reductive
groups over a finite field\bookinfo Ann.Math.Studies 107
\publ Princeton U.Press\yr1984\endref
\ref\key{L86}\by G.Lusztig\paper Sur les cellules gauches
des groupes de Weyl\jour C.R.Acad.Sci.Paris(A)\vol302
\yr1986\pages5-8\endref
\ref\key{L15}\by G.Lusztig\paper On conjugacy classes in a
reductive group\inbook Representations of Reductive Groups
\bookinfo Progr.in Math. 312\publ Birkh\"auser\yr2015
\pages333-363\endref
\ref\key{L19}\by G.Lusztig\paper
A new basis for the representation ring of a Weyl group\jour
Repres.Th.\vol23\yr2019\pages439-461\endref
\ref\key{L20}\by G.Lusztig\paper The Grothendieck group of
unipotent representations: a new basis\jour
Represent.Th.\vol24\yr2020\pages178-209\endref
\ref\key{L23}\by G.Lusztig\paper Precuspidal families and
indexing of Weyl group representations\lb\jour
arxiv:2304.05130\endref
\ref\key{LS}\by G.Lusztig and E.Sommers\paper
Constructible representations and Catalan numbers           jour arxiv:2403.02550\endref
\endRefs
\enddocument